\documentclass[reqno,12pt]{amsart}

\title[Phase unwinding, decompositions of Hardy spaces]{Phase
  unwinding, or invariant subspace decompositions of Hardy spaces}

\author{Ronald R.~Coifman}\address{Department of Mathematics, Program
  in Applied Mathematics, Yale University, New Haven, CT 06510, USA}
\email{coifman-ronald@yale.edu}

\author{Jacques Peyri\`ere} \address{Laboratoire de Math\'ematiques
  d'Orsay, Univ.~Paris-Sud, CNRS, Universit\'e Paris-Saclay, 91405
  Orsay, France.} \email{jacques.peyriere@math.u-psud.fr}

\keywords{Blaschke factorization, Phase unwinding, Takenaka basis,
  Hardy spaces, inner function, invariant subspaces, multiscale
  decomposition.}

\subjclass[2000]{30B50; 30A10, 42C40, 65T99}
\usepackage{amsmath,amsfonts,amssymb,mathrsfs} 
\newtheorem{theorem}{Theorem}
\newtheorem{lemma}[theorem]{Lemma}
\newtheorem{corollary}[theorem]{Corollary}
\theoremstyle{definition}
\theoremstyle{remark}
\newtheorem*{remark*}{Remark}
\newcommand{\abs}[1]{\lvert #1\rvert}
\newcommand{\dif}{\mathrm{d}}

\newcommand{\e}{\mathrm{e}}
\newcommand{\mi}{\mathrm{i}}

\newcommand{\sca}[1]{\left\langle #1 \right\rangle }
\newcommand{\cnj}[1]{\overline{#1}}
\newcommand{\hardy}[1]{\mathsf{H}^{#1}}
\newcommand{\hardyT}[1]{\mathsf{H}^{#1}({\mathbb T})}
\newcommand{\hardyR}[1]{\mathsf{H}^{#1}({\mathbb R})}
\newcommand{\uhp}{{\mathbb H}}
\newcommand{\InvSp}{{\mathsf V}}

\newcommand{\scale}{{\Delta}}
\begin{document}

\begin{abstract}
We consider orthogonal decompositions of invariant subspaces of Hardy
spaces, these relate to the Blaschke based phase unwinding
decompositions ~\cite{coifman,CSW,nahon}.

We prove convergence in $L^p$. In particular we build an explicit
multiscale wavelet basis. We also obtain an explicit unwindinig
decomposition for the singular inner function, $ \exp 2\mi\pi /x$.

\end{abstract}

\maketitle
\section{Introduction}
 Our goal is to extend and clarify convergence properties of the phase
 unwinding expansions in ~\cite{nahon,coifman,CSW,qian} as well as
 expansions obtained by the algorithm of adaptative Fourier
 decomposition~\cite{qian1,qian2,qian3,qian} where each function in
 $\hardyT2$, admits its own adapted (unwound) decomposition in an
 orthonormal system of basis functions consisting of partial products
 of Blaschke products. We extend the result to $\hardyT{p}$ for $p\in
 (1,+\infty)$.  We also discuss the relation to various
 generalizations of the Takenaka Malmquist bases, both for the Torus
 and the upper half plane.  In particular we show that there is a
 multiscale analysis of $\hardy{2}({\mathbb R})$, and that, at each
 level, there is a function whose translates make an orthonormal
 basis. This in the same spirit as studies of hyperbolic
 wavelets~\cite{pap-schipp,pap,feichtinger}.  The main difference is
 that we use different grids, which allows to get a formalism very
 close to wavelets. More precisely, let
 \begin{equation*}
\phi(x) = \frac{\Gamma(x-1+\mi)}{\sqrt{\pi}\Gamma(x-\mi)}.
 \end{equation*}
 Then the functions $\phi(2^nx+j)\Delta(2^nx)$, for $n$ and $j$ in
 ${\mathbb Z}$, where~$\Delta$ is a suitable inner function, form an
 orthonormal basis of~$\hardy{2}$. Moreover, $\Delta(x) =
 \e^{2\mi\pi\bigl(x+\gamma(x)\bigr)}$, where $\gamma$ is an entire
 function, real on the real line, and such that $|\gamma(x)|\le
 0.004\,\pi |x|$.\medskip

 We also give an explicit unwinding of the singular inner function
 $\exp \frac{2\mi\pi}{x}$:
 \begin{equation*}
\exp \frac{2\mi\pi}{x} = \e^{-2\pi} + \bigl( 1-\e^{-4\pi}\bigr) \sum_{n\ge 0}
(-1)^n\e^{-2n\pi} B(x)^{n+1},
 \end{equation*}
 where $B$ is a Blaschke product whose zeros are $\{1/(j+\mi)\}_{j\in
   {\mathbb Z}}$.
 
\section{Preliminaries and notation}

For $p\ge 1$, $\hardyT{p}$ stands for the space of analytic functions~$f$
on the unit disk~${\mathbb D}$ such that
\begin{equation*}
\sup_{0< r<1} \int_{0}^{2\pi} |f(r\e^{\mi\theta})|^p\frac{\dif
  \theta}{2\pi} < +\infty.
\end{equation*}
Such functions have boundary values almost everywhere, and the Hardy
space $\hardyT{p}$ can be identified with the set of $L^p$ functions
on the torus~${\mathbb T}=\partial{\mathbb D}$ whose Fourier
coefficients of negative order vanish.

A subspace of~$\hardyT{p}$ is \emph{invariant} if it is invariant by
multiplication by~$\e^{\mi\theta}$ (or by~$z$, depending whether these
functions are considered as functions on~${\mathbb T}$ or~${\mathbb
  D}$). An inner function is a bounded analytic function on the unit
disk whose boundary values have modulus~1 almost everywhere. It is
known that the invariant subspaces are of the form $u\hardyT{p}$
where~$u$ is an inner function. The inner function~$u$ is determined
by the invariant subspace up to multiplication by a constant of
modulus~1.

Any $f\in \hardyT{p}$ decomposes as $gu$, where $u$ is inner and $g$
outer. The inner function in its turn can be
further decomposed as $BS$, where $B$ is a Blaschke product, which
accounts for all the zeros, and $S$ a singular inner
function~\cite{hoffman,helson}.

If $f$ and $g$ are two functions on ${\mathbb T}$ (in $L^p$ and
$L^{p/(p-1)}$ for some $p\in [1,+\infty)$), let
\begin{equation*}
\sca{f,g} = \frac{1}{2\pi} \int_0^{2\pi}
f(\e^{\mi\theta})\cnj{g(\e^{\mi\theta})} \,\dif \theta.
\end{equation*}

Let $H$ be the operator of orthogonal projection of $L^2({\mathbb T})$
onto $\hardyT2$. 
It results from the properties of the Hilbert that this operator
extends as a bounded operator from $L^p({\mathbb T})$ to $\hardyT{p}$
for $1< p< +\infty$.

If $u$ is an inner function, let $\chi_u$ be the operator of
multiplication by~$u$ (which is an isometry of all the $L^P$). Then the
operator $H_u=\chi_uH\chi_u^{-1}$ is the operator of orthogonal projection
of $L^2$ onto $u\hardyT2$. It results that this operator extends as a
bounded operator from $L^p({\mathbb T})$ to $\hardyT{p}$ for all $p\in
(1,+\infty)$ with a norm independent of~$u$. In other terms, for all~$p>
1$, there exists~$C_p$ such that, for all~$u$ and all $f\in
L^{p}({\mathbb T})$,
\begin{equation}\label{projection}
\|H_uf\|_p \le C_p\|f\|_p.
\end{equation}
\bigskip

There is a parallel theory for analytic functions on the upper half
plane $\uhp = \{x+\mi y \ :\ y>0\}$. The space of analytic
functions~$f$ on $\uhp$ such that
$$\sup_{y>0} \|f(\cdot+\mi y)\|_{L^p({\mathbb R})} < +\infty$$
is denoted by $\hardyR{p}$. These functions have boundary values in
$L^p({\mathbb R})$ when $p\ge 1$. The space $\hardyR{p}$ is identified
to the space of $L^p$ functions whose Fourier transform vanishes on
the negative half line~$(-\infty,0)$.

A subspace of $\hardyR{2}$ is said to be invariant if it is stable by
multiplication by the functions $\e^{2\mi\pi\xi x}$ for all
$\xi>0$. As previously, the invariant subspaces are of the form
$u\,\hardy{2}$ where $u$ is an inner function, i.e., an analytic
function on~$\uhp$ whose boundary values are of modulus~1 almost
everywhere.

As previously, the operators of orthogonal projections on invariant
subspaces extend, for any $p\in (1,+\infty)$, as continuous operators
on $\hardyR{p}$ with a uniform bound for their norms.

\section{Phase unwinding on the torus}\label{Unwinding}

In this section, one simply writes $\hardy{p}$ instead of $\hardyT{p}$.
\subsection{Phase unwinding}

The following construction is a slight generalization of the one
described in~\cite{coifman,CSW,nahon}.

One starts with $f\in \hardy{p}$. We choose a projector~$Q_0$ on some
subspace of~$\hardy{p}$ and write $g_0=Q_0f$ and $f = g_0+u_1f_1$,
where $u_1$ is inner and $f_1$ outer. Then choose a projector~$Q_1$,
not necessarily different from~$Q_0$, and write $g_1=Q_1f_1$ and
$f_1=g_1+u_2f_2$, and so on indefinitely unless reaching~0. This leads
to the expansion
\begin{equation}\label{unwinding}
f = g_0+u_1g_1+u_1u_2g_2+\cdots+g_{n-1}\prod_{1\le j<n}u_j +
f_n\prod_{1\le j\e n}u_j,
\end{equation}
which is orthogonal when~$p=2$.

When the process does not stop, it is natural to ask in what sense the
series
\begin{equation}\label{series}
g_0+\sum_{n\ge 1}g_n\!\prod_{0\le k\le n}u_k
\end{equation}
represents~$f$. We shall answer this question in the next section.
\bigskip

In previous works, the projectors~$Q_j$ were of the form
$\mathsf{Id}-H_v$, where $v$ is an inner function. Let us investigate
this case.

Let $(v_j)_{j\ge 1}$ be a sequence of inner functions.
Let $f$ in $\hardy{p}$ for some $p\ge 1$. Define by recursion three
sequences (maybe finite) of functions $(f_n)_{n\ge 0}$, $(g_n)_{n\ge
  0}$, and $(u_n)_{n\ge 1}$, where, for $n\ge 1$, the $u_n$ are inner,
and the $f_n$ outer:
\begin{itemize} 
\item[--] $f_0=f$,
\item[--] to pass from step~$n$ to step~$n+1$, consider the projection
  $h_n = H_{v_{n+1}}f_n$ of~$f_n$ on $v_{n+1} \hardy{p}$; if~$h_n=0$
  then stop, otherwise let $f_{n+1}$ be the outer part of~$h_n$ and
  $u_{n+1}$ its inner part, and set $g_n = f_n-h_n =
  f_n-u_{n+1}f_{n+1}$.
\end{itemize}
Notice that, since $h_n\in v_{n+1}\hardy{p}$, $u_{n+1}/v_{n+1}$ is an
inner function. In particular $u_{n+1}$ vanishes at the zeros
of~$v_{n+1}$ with a not smaller multiplicity.\medskip

When $v_{n+1}$ is a convergent Blaschke products, consider a
Malmquist-Takenaka basis $e_1,\,e_2,\cdots$ associated with this
product (see~\cite{takenaka} and the next section). Then, for $f\in
\hardy{2}$, $g_n = \sum_{j} \sca{f_n,e_j}e_j$. As the functions $e_j$
are bounded, these scalar products are well defined if $f\in
\hardy{p}$ and the above expression for~$g_n$ holds as well.\medskip

The above construction, when it does not stop, appeals further
comments. Any zero of~$v_n$ is also a zero of~$u_n$ (with a
multiplicity not smaller). Consider the
following decreasing sequence of subspaces of $\hardy{p}$:
$$\InvSp_n = u_1u_2\cdots u_n\hardy{p}$$
and the space $\InvSp_{\infty} = \displaystyle \bigcap_{n\ge
  0}\InvSp_n$.

If a function $h$ is in $\InvSp_{\infty}$, the set of its zeros
contains all the zeros of the $v_n$ counted with their multiplicities
at least. This means that if $\sum |1-z_j|=\infty$ (where the $z_j$
are the zeros of the $(v_n)_{n\ge 1}$ repeated according to their
multiplicities) $h=0$. In other word, under this hypothesis,
$\InvSp_{\infty} = \{0\}$.
\medskip

Now, let us describe a few choices of the $v_j$.
If we take $v_j(z)=z$ for all z, then $g_n(z) = f_n(0)$. This is a
variant of the case studied in~\cite{nahon,coifman,CSW}. Indeed in
these articles, in the recursion $f_n-f_n(0)$ is decomposed as gSB,
where $g$ is outer, $S$ singular inner, and $B$ a
Blaschke product; then one sets $f_{n+1} =gS$ instead of $f_{n+1}=g$.

Another possibility is to take $v_n(z) =
\frac{z-a_n}{1-\overline{a}_nz}$ with $|a_n|< 1$ and $\sum (1-|a_n|) =
\infty.$ Then $g_n(z) = f_n(a_n)(1-|a_n|^2)/(1-\overline{a}_nz)$.

An important remark is that the sequence $v_j$ can be defined 'on
the fly' to be well adapted to the function~$f$ under analysis. For
instance, one can choose $v_{n+1}$ within a collection of Blaschke
products to maximize the $L^2$-norm of~$g_n$. This is the greedy
algorithm of~Qian et al.~\cite{qian0,qian1,qian2}.

Here is a variant of the above construction. Let $a_j$ be a sequence
of numbers of moduli less than~1. Given a~$f$ one considers the
recursion:
\begin{itemize}
\item[--] $f_0=f$
\item[--] for $n\ge 1$, $c_{n-1}=f_{n-1}(z_{n-1})$, $f_{n-1}- c_{n-1} =
  u_nf_n$, where $u_n$ is inner and~$f_n$ outer.
\end{itemize}

\subsection{Nested invariant subspaces}

\begin{theorem}\label{convergence}
 Let $(\InvSp_n)_{n\ge 0}$, with $\InvSp_0 = \hardy{2}$, be a
 decreasing sequence of invariant subspaces. Set $\InvSp_{\infty} =
   \bigcap \InvSp_n$ and let ${\mathcal P}_n$ stand for the operator
   associated with the inner function defining $\InvSp_n$. Then, for
   all~$p\in (1,+\infty)$ and $f\in \hardy{p}$, one has
$$\lim_{n\to +\infty}\|{\mathcal P}_nf-{\mathcal P}_{\infty}f\|_p =
  0.$$
\end{theorem}

\proof Fix $p\in (1,+\infty)$ ($p\ne2$) and $p_0\in (1,+\infty)$ such
that $p$ lies in the open interval delimited by $2$ and $p_0$. Let
$g\in \hardy{2}\cap \hardy{p_0}$. By H\"older inequality, there exists
$\alpha\in (0,1)$, depending only on~$p$ and~$p_0$, such that
$\|{\mathcal P}_ng-{\mathcal P}_{\infty}g\|_p \le \|{\mathcal
  P}_ng-{\mathcal P}_{\infty}g\|_2^\alpha\|{\mathcal P}_ng-{\mathcal
  P}_{\infty}g\|_{p_0}^{1-\alpha}$. It results from~\eqref{projection}
that
$$\|{\mathcal P}_ng-{\mathcal P}_{\infty}g\|_p \le
(2C_p)^{1-\alpha}\|{\mathcal P}_ng-{\mathcal
  P}_{\infty}g\|_2^\alpha\|g\|_{p_0}^{1-\alpha}$$
and
$$\lim_{n\to
  +\infty}\|{\mathcal P}_ng-{\mathcal P}_{\infty}g\|_p= 0.$$

Now, if $f\in \hardy{p}$, for all $g\in \hardy{2}\cap \hardy{p_0}$, one
has

\begin{equation*}
\|{\mathcal P}_nf-{\mathcal P}_{\infty}f\|_p \le \|{\mathcal P}_ng -
  {\mathcal P}_{\infty}g\|_p + \|{\mathcal P}_n(f-g)-{\mathcal
    P}_{\infty}(f-g)\|_p,
\end{equation*}
therefore (due to~\eqref{projection})
\begin{equation*}
\limsup_{n\to \infty}\|{\mathcal P}_nf-{\mathcal P}_{\infty}f\|_p \le
2C_p\inf_{g\in \hardy{2}\cap \hardy{p_0}} \|f-g\|_p = 0.
\end{equation*}

\begin{corollary}
Let ${\mathcal Q}_n={\mathcal P}_n-{\mathcal P}_{n+1}$. Then, for
all~$p\in (1,+\infty)$ and $f\in \hardy{p}$, the series
$$\sum_{n\ge 0}{\mathcal Q}_nf$$
converges to $f-{\mathcal P}_{\infty}f$ in $L^p$.
\end{corollary}

In particular, this proves that the series~\eqref{series} converges
to~$f$ in $\hardy{p}$ provided that $1< p< +\infty$.\medskip

This corollary also contains the theorem by Szab\'o~\cite{szabo} and
by Qian et al.~\cite{qian3} on the $H^p$-convergence of
Malmquist-Takenaka series.

\subsubsection{Malmquist-Takenaka bases}

For the reader's convenience we give an account of Malmquist-Takenaka
bases.

\begin{lemma}\label{takenaka}
Let $a$ be a complex number of modulus less than~1. Then
$(z-a)\hardy{2}$ has codimension~1 in $\hardy{2}$ and $\displaystyle
\frac{\sqrt{1-|a|^2}}{1-\cnj{a}z}$ is a unit vector in the orthogonal
complement of $(z-a)\hardy{2}$ in $\hardy{2}$.
\end{lemma}

\proof One has
\begin{eqnarray*}
\sca{(z-a)f(z),\frac{\sqrt{1-|a|^2}}{1-\cnj{a}z}} &=& \frac{1}{2\pi}
\int_{-\pi}^{\pi} (\e^{\mi\theta}-a)f(\e^{\mi\theta})
\frac{\sqrt{1-|a|^2}}{1-a\e^{-\mi\theta}}\,\dif \theta\\
&=& \frac{\sqrt{1-|a|^2}}{2\pi}\int_{-\pi}^{\pi}
\e^{\mi\theta}f(\e^{\mi\theta})\,\dif \theta = 0. 
\end{eqnarray*}

Also, if $f$ is orthogonal to $(1-\cnj{a}z)^{-1}$ one has
\begin{equation*}
0 = \frac{1}{2\pi} \int_{-\pi}^{\pi}
\frac{f(\e^{\mi\theta})}{1-a\e^{-\mi\theta}}\,\dif \theta\\
= \frac{1}{2\mi\pi} \oint \frac{f(z)}{z-a}\,\dif z = f(a),
\end{equation*}
so $f\in (z-a)\hardy{2}$.\medskip

Now $(a_n)_{n>0}$ is a sequence of complex numbers of moduli less
than~1 such that
\begin{equation}\label{convergence}
  \sum_{n\ge 1}(1-|a_j|^2) = +\infty.
\end{equation}
For $n\ge 0$, let
 $$B_n(z) = \prod_{0\le j< n} \frac{z-a_j}{1-\cnj{a}_jz} \text{\quad
  and\quad} \phi_n(z) = B_n(z)\frac{\sqrt{1-|a_n|^2}}{1-\cnj{a}_nz}.$$
It results from Lemma~\ref{takenaka} that the functions $\phi_n$ form
an orthonormal basis of $\hardy{2}$. Indeed, the spaces
$B_n\hardy{2}$, for $n\ge 0$, form a nested sequence of invariant
subspaces, and $\phi_n$ is a basis of the unidimensional space
$B_n\hardy{2}\ominus B_{n+1}\hardy{2}$. The bases so obtained are the
\emph{Malmquist-Takenaka bases}. Theorem~\ref{convergence} implies
that, if $1< p<+\infty$ and $f\in \hardy{p}$, the series
$\displaystyle \sum_{n\ge 0} \sca{f,\phi_n}\phi_n$ converges to~$f$ in
$\hardy{p}$.

\section{The upper half plane}

In this section, one simply writes $\hardy{p}$ instead of
$\hardyR{p}$.

\subsection{Malmquist-Takenaka bases}

Among the inner functions~$u$ there are the Blaschke products: let
$(a_j)_{1\le j}$ be a sequence (finite or not) of complex numbers
with positive imaginary parts and such that
\begin{equation}\label{upper} 
\sum_{j\ge 0} \frac{\Im a_j}{1+|a_j|^2} < +\infty.
\end{equation}
The corresponding Blaschke product is
$$B(x) = \prod_{j\ge 0}
\frac{\abs{1+a_j^2}}{1+a_j^2}\,\frac{x-a_j}{x-\cnj{a}_j},$$
where, $0/0$, which appears if $a_j=\mi$, should be understood as~1.
The factors $\displaystyle \frac{\abs{1+a_j^2}}{1+a_j^2}$ insure the
convergence of this product when there are infinitely many
zeroes. But, in some situations, it is more convenient to use other
convergence factors as we shall see below.

Whatever the series~\eqref{upper} be convergent or not, one defines
(for $n\ge 0$) the functions
\begin{equation*}
\phi_n(x) = \frac{1}{\sqrt{\pi}}\left( \prod_{0\le j< n}
\frac{x-a_j}{x-\cnj{a}_j}\right)\, \frac{1}{x-\cnj{a}_n}.
\end{equation*}
Then these functions form a orthonormal system in $\hardy{2}$. If the
series~\eqref{upper} diverges, it is a basis of $\hardy{2}$, otherwise
it is a basis of the orthogonal complement of $B\,\hardy{2}$ in
$\hardy{2}$.

For $1< p< +\infty$, and $f\in \hardy{p}$, the series $\sum_{n\ge 0}
\sca{f,\phi_n}\phi_n$ converges in $\hardy{p}$ (towards~$f$ if the
series~\eqref{upper} diverges). The proof is the same as previously.

\subsection{A multiscale decomposition}

The infinite products
\begin{equation}
G_n(x) = \prod_{j\le n} \frac{j-\mi}{j+\mi}\, \frac{x-j-\mi}{x-j+\mi}
\text{\quad and\quad } G(x) = \prod_{j\in {\mathbb Z}} \frac{j-\mi}{j+\mi}\,
\frac{x-j-\mi}{x-j+\mi}
\end{equation}
converge. As $\displaystyle
\frac{j-\mi}{j+\mi}\times\frac{-j-\mi}{-j+\mi} = 1$, one has
\begin{equation*}
G(x) = -\lim_{n\to +\infty} \prod_{|j|\le n} \frac{x-j-\mi}{x-j+\mi},
\end{equation*}
which shows that $G$ is periodic of period~1. It appears that these
products can be expressed in terms of known functions.

\begin{lemma} We have
\begin{equation*}
G_n(x) = \frac{\Gamma(-\mi-n)}{\Gamma(\mi-n)}\,
\frac{\Gamma(x-n+\mi)}{\Gamma(x-n-\mi)} \text{\quad and\quad } G(x) =
\frac{\sin \pi(\mi-x)}{\sin \pi(\mi+x)}.
\end{equation*}
\end{lemma}

\proof The well known formula $\displaystyle \Gamma(z) = \lim_{n\to
  +\infty} \frac{n!n^z}{z(z+1)\cdots(z+n)}$ yields the expression of
$G_0$.

On the other hand,
\begin{eqnarray*}
\frac{1}{G_0(-x)} &=& \prod_{j\le 0}
\frac{j+\mi}{j-\mi}\,\frac{-x-j+\mi}{-x-j-\mi} = \prod_{j\ge 0}
\frac{-j+\mi}{-j-\mi}\,\frac{-x+j+\mi}{-x+j-\mi}\\
&=& \prod_{j\ge 0} \frac{j-\mi}{j+\mi}\,\frac{x-j-\mi}{x-j+\mi} =
-\frac{x-\mi}{x+\mi} \prod_{j\ge 1}
\frac{j-\mi}{j+\mi}\,\frac{x-j-\mi}{x-j+\mi}.
\end{eqnarray*}

Therefore
\begin{eqnarray*}
G(x) &=& -\frac{(x+\mi)G_0(x)}{(x-\mi)G_0(-x)} =
-\frac{-(x+\mi)\Gamma(x+\mi)\Gamma(-x-\mi)}
{(-x+\mi)\Gamma(x-\mi)\Gamma(-x+i)}\\
&=& -\frac{\Gamma(x+\mi)\Gamma\bigl(1-(x+\mi)\bigr)}
{\Gamma(x-\mi)\Gamma\bigl(1-(x-\mi)\bigr)} = -\frac{\sin
  \pi(x-\mi)}{\sin \pi(x+\mi)}.
  \end{eqnarray*}

\subsection{An orthonormal system}

Consider the function
$$\phi(x) = \frac{\Gamma(x-1+\mi)}{\sqrt{\pi}\Gamma(x-\mi)}.$$
It is easily checked that
$$\phi(x-n) =
\frac{\Gamma(\mi-n)}{\Gamma(-\mi-n)} \, \frac{G_n(x)}{\sqrt{\pi}\bigl(
  x-(n+1)+\mi\bigr)}.$$
Set $\phi_n(x) = \phi(x-n)$. For fixed~$m$, the functions
$\phi_n/G_m$, for $n\ge m$, form a Malmquist-Takenaka basis of
$(G/G_m)\hardy{2}$. In other terms, the functions $\phi_n$, for $n\ge
m$, form an orthonormal basis of $G_m\hardy{2}\ominus
G\hardy{2}$. This means that the functions $\phi(x-n)$ (for $n\in
{\mathbb Z}$) form a Malmquist-Takenaka basis of the orthogonal
complement of $G\hardy{2}$ in $\hardy{2}$.

Then the same proof as the one of Theorem~\ref{convergence} yields the
following results.

\begin{lemma} Let $f\in \hardy{p}$ for some $p> 1$. Then both series
$$\sum_{j<0} \sca{f,\phi_j}\phi_j \text{\quad and\quad} \sum_{j\ge 0}
  \sca{f,\phi_j}\phi_j$$
are convergent in $\hardy{p}$ and
$$f = \sum_{n\in {\mathbb Z}} \sca{f,\phi_n}\phi_n + {\mathcal
  P}_pf,$$
where ${\mathcal P}_p$ stands for the extension to $L^p$ of the
orthogonal projector on $G\hspace{1pt}H^2$.
\end{lemma}

\subsubsection{Multiscale decomposition}

As $|1-G(x)|\le C\min \{1,|x|\}$ the product
\begin{equation}
\scale(x) = \prod_{j< 0} G(2^jx)
\end{equation}
is convergent and $\displaystyle \lim_{n\to -\infty} \scale(2^nx) = 1$
uniformly on compact sets.


Consider the following subspaces of $\hardy{2}$:
\begin{equation*}
\InvSp_n = \scale(2^nx)\hardy{2}.
\end{equation*}
This is a decreasing sequence. The space $\displaystyle
\InvSp_{+\infty} = \bigcap_{n\in {\mathbb Z}} \InvSp_n$ is equal
to~$\{0\}$ since a non-zero function in this space would have too many
zeros, and the space $\displaystyle \InvSp_{-\infty} =
\mathrm{closure~of}\bigcup_{n\in {\mathbb Z}} \InvSp_n$ is equal to
$\hardy{2}$ since $\scale(2^nx)$ converges to 1 uniformly on compact
sets when $n$
goes to $-\infty$.

For all $n$ and~$j$, let
\begin{equation}\label{wavelets}
  \phi_{n,j}(x) = 2^{n/2}\phi(2^nx-j)\scale(2^{n}x).
\end{equation}
Then, for all~$n$, $(\phi_{n,j})_{j\in {\mathbb Z}}$ is an orthonormal
basis of $\InvSp_{n}\ominus \InvSp_{n+1}$. At last
$(\phi_{n,j})_{n,j\in{\mathbb Z}}$ is an orthonormal basis of
$\hardy{2}$. The following theorem results from the preceding
discussion.

\begin{theorem}
Let $\prec$ stand for the lexicographic order on ${\mathbb
  Z}\times{\mathbb Z}$. Then, if $f\in \hardy{p}$ for some $p> 1$,
one has
$$\lim_{\substack{(n,j)\to \infty\\\hardy{p}}} \sum_{(m,k)\prec (n,j)}
\sca{f,\phi_{m,k}}\phi_{m,k} = f.$$
\end{theorem}
\bigskip

Let us give another expression of $\scale$. Write
$\scale=\e^{\mi\psi}$. A simple calculation yields
$$G(x) = \exp 2\mi\left( \pi x + \sum_{n\ge 1} \frac{\e^{-2\pi n}}{n}
\sin \pi nx\right).$$
Then
\begin{eqnarray*}
\psi(x) = 2\pi x + 2\sum_{n\ge 1}\frac{e^{-2\pi n}}{n} \sum_{k\ge 1}
\sin 2\pi2^{-k}nx.
\end{eqnarray*}
So, if we set
$$\Xi(t) = \sum_{k\ge 1} \sin 2^{-k}t$$
we get
$$ \psi(x) = 2\pi x + 2\sum_{n\ge 1} \frac{\e^{-2\pi n}}n\Xi(2\pi nx).$$
\medskip

It is worth noticing that $\psi$, as a function on ${\mathbb R}$, is
increasing and that $|\psi(x)-2\pi x| <
2\pi|x|\e^{-2\pi}/(1-\e^{-2\pi})<0.004\,\pi|x|$. This means that
$\scale$ is ``not far'' from being periodic of period~1, and that
\eqref{wavelets} is reminiscent of the usual formula for wavelets.

A more precise bound of the form\\ $|\psi(x)-2\pi x|\le
c\min\{|x|,\log (1+|x|)\}$ can be obtained: one has $ |\Xi(x)|\le n+2$
when $n\ge 0$ and $2^n\le |x|\le 2^{n+1}$.

\section{Explicit phase unwindings}

The unwinding procedure can be performed for $\hardy{2}({\mathbb R})$
functions in the same way as in Section~\ref{Unwinding}. This time, we
take $h_j(x) = \displaystyle \frac{1}{\sqrt{\pi}(z-\overline{z}_j)}$, with
$\Re z_j> 0$.\medskip

Let $f(x) = \e^{2\mi\pi x}$. The functions $f-f(\mi)$ and $G$ have the
same zeros. So we can write $f(x)-f(\mi) = G(x)h(x)$:
\begin{eqnarray*}
h(x) &=& \frac{\e^{-2\pi}\bigl(\e^{2\mi\pi(x-i)}-1\bigr)\sin
  \pi(\mi+x)}{\sin \pi(\mi-x)}\\
&=& -2\mi\, \e^{\mi\pi(x-i+2\mi)} \sin \pi(\mi+x)\\
&=& 1-\e^{2\mi\pi(x+\mi)} = 1-\e^{-2\pi}f(x).
\end{eqnarray*}
It results the following identity
\begin{equation}\label{recur}
f = \frac{\e^{-2\pi}+G}{1+\e^{-2\pi}G}.
\end{equation}
So we get the following representation of~$\e^{2\mi\pi x}$,
\begin{equation}\label{ProUnwinding}
\e^{2\mi\pi x} = \e^{-2\pi} + \bigl( 1-\e^{-4\pi}\bigr) \sum_{n\ge 0}
(-1)^n\e^{-2n\pi} G(x)^{n+1},
\end{equation}
which gives an explicit unwinding series for $\displaystyle
\frac{\e^{2\mi\pi x}}{\sqrt{\pi}(x+\mi)}$:
\begin{equation}\label{CompleteUnwinding}
\frac{\e^{2\mi\pi x}}{\sqrt{\pi}(x+\mi)} =
\frac{\e^{-2\pi}}{\sqrt{\pi}(x+\mi)} + 
\bigl( 1-\e^{-4\pi}\bigr) \sum_{n\ge 0} (-1)^n\e^{-2n\pi}
\frac{G(x)^{n+1}}{\sqrt{\pi}(x+\mi)}.
\end{equation}

Also, by replacing $x$ by
$\frac{\mi}{2\pi}\frac{1+\e^{\mi\theta}}{1-\e^{\mi\theta}}$ in
  Formula~\eqref{ProUnwinding} we get an unwinding series of the singular
  inner function $\displaystyle \exp
  \frac{1+\e^{\mi\theta}}{1-\e^{\mi\theta}}$ on the torus.\medskip

In the same way, more general formulas can be obtained. Indeed, for
$\alpha>0$ let us define the following Blaschke product
$$G_\alpha(x) = \prod_{n\in {\mathbb Z}}
\frac{n-\mi\alpha}{n+\mi\alpha}\,
\frac{x-n-\mi\alpha}{x-n+\mi\alpha}.$$

We then have $G_{\alpha}(x) = \displaystyle \frac{\sin
  \pi(\mi\alpha-x)}{\sin \pi(\mi\alpha+x)}$ and
$$\e^{2\mi\pi x} =
\e^{-2\pi\alpha}+G_\alpha(x)\bigl(1-\e^{-2\pi\alpha}\e^{2\mi\pi
  x}\bigr).$$

This last formula leads to a variety of unwinding series: let
$(\alpha_n)_{n>0}$ be a sequence of positive numbers, set $a_n =
\e^{-2\mi\alpha_n}$; then
\begin{multline*}
\e^{2\mi\pi x} = a_1+(1-a_1a_2)G_{a_1}(x) -
a_1(1-a_2a_3)G_{a_1}(x)G_{a_2}(x) \\
+a_1a_2(1-a_3a_4)G_{a_1}(x)G_{a_2}(x)G_{a_3}(x)\\
-a_1a_2a_3(1-a_4a_5)G_{a_1}(x)G_{a_2}(x)G_{a_3}(x)G_{a_4}(x) + \cdots
\end{multline*}

Also, one can use identity~\eqref{recur} to analyse the singular inner
function on the real line defined by the Dirac mass at~0, namely
$\e^{-2\mi\pi/x}$:
\begin{equation*}
\e^{-\frac{2\mi\pi}{x}} = \frac{\e^{-2\pi}+B(x)}{1+\e^{-2\pi}B(x)},
\end{equation*}
where $B$ is the Blaschke product $G(1/x)$ whose zeros are the points
$-1/(n+\mi), \quad n\in{\mathbb Z}$.

\end{document}